\documentclass{article}
\usepackage{amssymb}
\usepackage{amsfonts}
\usepackage{amsmath}

\newtheorem{theorem}{Theorem}

\newtheorem{corollary}[theorem]{Corollary}

\newtheorem{definition}[theorem]{Definition}

\newtheorem{lemma}[theorem]{Lemma}

\newtheorem{proposition}[theorem]{Proposition}

\input{tcilatex}
\begin{document}

\title{Nilpotency and higher order derivatives\\
in differential geometry}
\author{Erc\"{u}ment H. Orta\c{c}gil}
\maketitle

\begin{abstract}
We recall the relation between nilpotency and higher order derivatives
remarked in [O1] and give some simple estimates for the order of a Klein
geometry.
\end{abstract}

\section{Motivation and preliminaries}

A "higher order geometry" on a smooth manifold $M$ is a "curved" homogeneous
space $G/H=M$ of "high order". Clearly, the two concepts of "curvature" and
"order" are in need of clarification in this definition. We will first take
a look at the second concept, refering to [O1] for more details. Suppose a
Lie transformation group $G$ acts \textit{effectively }and \textit{%
transitively }on $M.$ We fix some $a,b\in M$ and choose some $g\in G$ which
maps $a$ to $b,$ i.e., $g(a)=b.$ For $k\geq 0,$ we define

\begin{equation}
H_{k}(a,b,g)\overset{def}{=}\{h\in G\mid j_{k}(h)^{a,b}=j_{k}(g)^{a,b}\}
\end{equation}%
where $j_{k}(h)^{a,b}$ denotes the $k$-jet of $h$ with source at $a$ and
target at $b=h(a),$ i.e., $j_{k}(h)^{a,b}$ is the Taylor expansion of $h$ up
to order $k$ at $a.$ In words, $H_{k}(a,b,g)$ consists of those
transformations of $G$ which agree with $g$ up to order $k$ at $a.$ Clearly, 
$H_{0}(a,a,g)=$ the stabilizer of $G$ at $a$ which we denote by $G_{a}.$
Obviously we have the inclusions of sets

\begin{equation}
\{g\}\subset ...\subset H_{k}(a,b,g)\subset ...H_{1}(a,b,g)\subset
H_{0}(a,b,g)
\end{equation}

In other words, the above sets approximate $g$ better and better with
increasing $k.$ The descending chain in (1) stabilizes at $k=0,$ i.e., $%
H_{k}(a,b,g)=H_{0}(a,b,g)$ for all $k\geq 0$ if $G$ acts simply transitively
on $M.$ Now we have the following

\begin{proposition}
If $G$ is connected and acts effectively, then there exists an integer $m$
satisfying

1) $H_{m}(a,b,g)=\{g\}$

2) $m$ does not depend on $a,b$ and $g.$

\begin{definition}
The smallest integer $m$ satisfying 1) of Proposition 1 is the order of $%
(G,M)$ and is denoted by $o(G,M).$
\end{definition}
\end{proposition}

In short, if $o(G,M)=m,$ then \textit{any} transformation of $G$ is \textit{%
globally determined on }$M$ once we give its derivatives up to order $m$ at 
\textit{any} point $a\in M$ and $m$ is the smallest integer with this
property. The situation is anologous to the one where an analytic function
is determined globally on its domain of definition once we give its
derivatives of \textit{all orders }(hence its power series) at some point.
In our case, finitely many derivatives suffice since $G$ is a finite
dimensional Lie group and nontrivial elements of $G$ "move some points",
i.e., $G$ can be "detected from its action" by effectiveness. Now $o(G,M)=0$
if and only if $G$ acts simply transitively on $M.$ At this point, we are
naturally led to ask the following question.

\begin{equation*}
\text{\textbf{Q}}\mathbf{\ }\text{: Give examples of actions }(G,M)\text{
with "large" }o(G,M).
\end{equation*}

The reader should not have much difficulty in giving examples where $%
o(G,M)=1,2$ but after some effort he/she will be soon convinced that it is a
nontrivial problem to give examples with $o(G,M)\geq 3$ and may even suspect
that such examples do not exist. To make progress with \textbf{Q}, we
linearize the problem and reduce it to pure algebra. We first recall that a
tangent vector $\sigma \in T_{e}(G)$ at the identity $e\in G$ defines a
vector field $\overline{\sigma }$ on $M$ as follows: Let $c(t)$ be a curve
passing through $e$ at $t=0$ with the tangent vector $\sigma $ at $e.$ Now
for $x\in M,$ $\overline{\sigma }(x)$ is defined as the tangent of the curve 
$c(t)\star x$ at $t=0$ where $\star $ denotes the action of $G.$ This
construction gives a map

\begin{eqnarray}
T_{e}(G) &\longrightarrow &\mathfrak{X}(M) \\
\sigma &\longrightarrow &\overline{\sigma }  \notag
\end{eqnarray}%
where $\mathfrak{X}(M)$ denotes the (infinite dimensional) Lie algebra of 
\textit{all }vector fields on $M.$ We observe that (3) is injective since
the action is effective. Therefore we can identify the Lie algebra $T_{e}(G)=%
\mathfrak{g}$ with its image $\overline{\mathfrak{g}}\subset $ $\mathfrak{X}%
(M).$ Note that the modern and old interpretations merge in (3): The modern
interpretation of $T_{e}(G)=\mathfrak{g}$ as the Lie algebra of the abstract
Lie group $G$ on the left hand side of (3) and the old interpretation of $%
\overline{\mathfrak{g}}\subset $ $\mathfrak{X}(M)$ as the infinitesimal
generators of the transformation group $G$ of $M$ on the right hand side of
(3). It was Lie who discovered that if $G$ is "finite dimensional", then $%
\overline{\mathfrak{g}}$ has an algebraic structure now called a finite
dimensional Lie algebra.

Now (3) has another important property: Obviously, we have the inclusion $\{%
\overline{\sigma }(x)\mid \sigma \in \mathfrak{g}\}\subset T_{x}(M)$ for any 
$x\in M.$ However, since the action is transitive, we have in fact equality,
i.e.,

\begin{equation}
T_{x}(M)=\{\overline{\sigma }(x)\mid \sigma \in \mathfrak{g}\}
\end{equation}

Now we fix $x=a\in M$ arbitrarly and define

$\bigskip $%
\begin{equation}
\overline{\mathfrak{g}}_{0}\overset{def}{=}\{\overline{\sigma }\mid 
\overline{\sigma }(a)=0\}\subset \overline{\mathfrak{g}}
\end{equation}%
We easily check that the subalgebra $\overline{\mathfrak{g}}_{0}\subset 
\overline{\mathfrak{g}}$ is the Lie algebra of the stabilizer $G_{a}.$ We
observe that $\overline{\sigma }(a)=$ the value of $\overline{\sigma }$ at $%
a $ $=0$-jet of $\overline{\sigma }$ at $a=j_{0}(\overline{\sigma })^{a}.$
For $k\geq 0,$ we now define

\begin{equation}
\overline{\mathfrak{g}}_{k}\overset{def}{=}\{\overline{\sigma }\mid j_{k}(%
\overline{\sigma })^{a}=0\}\subset \overline{\mathfrak{g}}
\end{equation}%
and obtain the filtration

\begin{equation}
\{0\}\subset ...\subset \overline{\mathfrak{g}}_{k}\subset ...\subset 
\overline{\mathfrak{g}}_{1}\subset \overline{\mathfrak{g}}_{0}\subset 
\overline{\mathfrak{g}}\subset \mathfrak{X}(M)
\end{equation}%
\ 

We observe that we differentiate smooth vector fields on $M$ to obtain the
filtration (7), i.e., \textit{we use calculus}. We can now prove

\begin{proposition}
The exists an integer $r$ satisfying $\overline{\mathfrak{g}}_{r}=0$ and the
smallest such integer $r$ satisfies either $r=o(G,M)$ or $r=o(G,M)+1.$
\end{proposition}

The problem of which case occurs in Proposition 3 depends on the behaviour
of some descrete subgroups. Since our purpose is to linearize \textbf{Q }and
therefore avoid global issues, henceforth we will concentrate on the
smallest integer $r$ given by Proposition 3. The geometric meaning of $r$ is
clear from its construction: Any vector field (= infinitesimal generator) in 
$\overline{\mathfrak{g}}$ is globally determined on $M$ by its $r$-jet at
any point $a\in M.$

Now a very crucial step: Using the identification of $\mathfrak{g}$ and $%
\overline{\mathfrak{g}}$ given by (3), we can pull back the
calculus-filtration (7) to a filtration inside $\mathfrak{g}=T_{e}(G)$ as

\begin{equation}
\{0\}=\mathfrak{g}_{r}\subset \mathfrak{g}_{r-1}\subset ...\subset \mathfrak{%
g}_{1}\subset \mathfrak{g}_{0}\subset \mathfrak{g}
\end{equation}%
where $\mathfrak{g}_{k}$ is by definition the preimage of $\overline{%
\mathfrak{g}}_{k}$ by (3). The question is now whether we can define (8) 
\textit{purely algebraically} avoiding (7)!!$\mathfrak{\ }$Remarkably, the
answer turns out to be affirmative and the construction goes as follows. Let
us interpret $\mathfrak{g}$ in (8) as an \textit{abstract }finite
dimensional Lie algebra, forgetting that it is the preimage of the
infinitesimal generators $\overline{\mathfrak{g}}$ in (7) arising from the
action $(G,M).$ Let $\mathfrak{g}_{0}\subset \mathfrak{g}$ be a subalgebra,
called the stabilizer subalgebra. We assume that the pair $(\mathfrak{g,g}%
_{0})$ (called an infinitesimal Klein geometry in [O1]) is effective, i.e., $%
\mathfrak{g}_{0}$ does not contain nontrivial ideals of $\mathfrak{g}$ (note
that since $G$ acts effectively on $M,$ the pair $(\overline{\mathfrak{g}},%
\overline{\mathfrak{g}}_{0})$ in (7) and therefore $(\mathfrak{g,g}_{0})$ in
(8) are effective). Starting with the abstract data $(\mathfrak{g,g}_{0}),$
we now define inductively

\begin{eqnarray}
\mathfrak{g}_{k+1}\overset{def}{=}\{x &\in &\mathfrak{g}_{k}\mid \lbrack x,%
\mathfrak{h]\subset g}_{k},\text{ \ }k\geq 0\} \\
&=&\{x\in \mathfrak{g}\mid \lbrack x,g\mathfrak{]\subset g}_{k},\text{ \ }%
k\geq 0\}  \notag
\end{eqnarray}

The second equality in (9) is easily seen by induction, one inclusion being
obvious. By definition $\mathfrak{g}_{k+1}\subset \mathfrak{g}_{k}\subset 
\mathfrak{g}$ and we easily check that the chain stabilizes at $\{0\}$ using
effectiveness, denoting the smallest such integer by $o(\mathfrak{g,g}_{0}).$
Note that we have an obvious representation of $\mathfrak{g}_{0}$ on $%
\mathfrak{g/g}_{0}$ and $\mathfrak{g}_{1}$ is the kernel of this
representation. Assuming $\mathfrak{g}_{k}$ is defined, $\mathfrak{g}_{k+1}$
is the kernel of the represantion of $\mathfrak{g}_{k}$ on $\mathfrak{g/g}%
_{k}.$ The filtration given by (9) is called the Weissfeiler filtration in
some works. To our knowledge, it appeared first in [GS]. Now we can prove

\begin{proposition}
The abstract filtration defined by (9) coincides with the
calculus-filtration (8).
\end{proposition}

In view of Proposition 4, \textbf{Q} becomes now a purely algebraic question
about Lie algebras:

\begin{equation*}
\text{\textbf{Q' : }How do we construct effective Lie algebra pairs }(%
\mathfrak{g,g}_{0})\text{ of large order?}
\end{equation*}

We called \textbf{Q' }"The first fundamental problem of higher order
geometry and jet theory (\textbf{FP1)}" in [O1]. Recalling that $o(\mathfrak{%
g,g}_{0})=0$ if and only if $\mathfrak{g}_{0}=\{0\},$ we now state

\begin{proposition}
Let $(\mathfrak{g,g}_{0})$ be an effective pair with stabilizer $\mathfrak{g}%
_{0}\neq \{0\}$ (so that $o(\mathfrak{g,g}_{0})\geq 1).$ If $\mathfrak{g}$
or $\mathfrak{g}_{0}$ is semisimple, then $o(\mathfrak{g,g}_{0})=1.$ The
same conclusion holds if one of $\mathfrak{g,g}_{0}$ is compact.
\end{proposition}

Anybody familiar with Lie algebras is surely familiar also with the
wonderful theory of semisimple algebras, their classification by root
systems, representations...etc. However, Proposition 5 warns us that the
pursuit of \textbf{Q' }will drive us in the \textit{opposite direction} of
the semisimple theory into the territory of nilpotent and solvable Lie
algebras.

Finally, a few words on curvature: A very subtle theory emerges if we give
"curvature" to a \textit{simply transitive} pair $(G,M)$ where the order is 
\textit{zero}. To give the reader an idea of the level of depth of this
theory, we will refer him/her to Part II of [O1] (see also [O2]). Therefore,
it is natural to expect that a "curved" $(G,M)$ of higher order will be an
even more subtle geometric structure. For instance, recall that a Lie group $%
G$ has a maximal compact subgroup $H$ and $G/H$ is topologically trivial.
This fact together with the second assertion of Proposition 5 suggests that
the curvature of a higher order structure contains information which can not
be detected by topological methods.

\section{Two simple bounds for order}

Our purpose in this section is to take a small step towards the solution of 
\textbf{Q'. }Henceforth, we will always assume that $(\mathfrak{g,g}_{0})$
is an effective pair with $\mathfrak{g}_{0}\neq 0$ so that $o(\mathfrak{g,g}%
_{0})\geq 1.$ In particular $\mathfrak{g}_{0}\subset \mathfrak{g}$ is a 
\textit{proper }subalgebra which is not an ideal. First, let us recall an
important property of the abstract filtration (9). It is easy to show that

\begin{equation}
\lbrack \mathfrak{g}_{i},\mathfrak{g}_{j}]\subset \mathfrak{g}_{i+j}\text{ \
\ \ \ }i,j\geq 0
\end{equation}

In particular, we have $[\mathfrak{g}_{0},\mathfrak{g}_{i}]\subset \mathfrak{%
g}_{i}$ so that $\mathfrak{g}_{i}$ is an ideal of not only $\mathfrak{g}%
_{i+1}$ but an ideal of $\mathfrak{g}_{0}.\mathfrak{\ }$Using (10) we easily
show that $\mathfrak{g}_{1}$ is nilpotent and therefore contained in the
largest nilpotent ideal (called nilradical) of $\mathfrak{g.}$

Now we recall $ad(x)(y)=[x,y]$ and for simplicity of notation, we denote the
linear map $ad(x_{k})\circ ...\circ ad(x_{1}):\mathfrak{g\longrightarrow g}$
by $ad(x_{k},...,x_{1})$ where we compose from right to left. Now we inspect
(9) more carefully. Suppose $r=o(\mathfrak{g,g}_{0})\geq 1$ so that $%
\mathfrak{g}_{r}=0$ but $\mathfrak{g}_{r-1}\neq 0.$ We choose some $a\in 
\mathfrak{g}_{r-1},$ $a\neq 0.$ If $[a,\mathfrak{g]\subset }\mathfrak{g}%
_{r-1},$ then by (8) we have $a\in $ $\mathfrak{g}_{r}=0,$ which is not
case. Therefore, there exists some $x_{1}\in \mathfrak{g}$ such that $%
[x_{1},a]\notin \mathfrak{g}_{r-1}.$ Since $[\mathfrak{g}_{0},\mathfrak{g}%
_{r-1}]\subset \mathfrak{g}_{r-1}$ by (10), we have $x_{1}\in \mathfrak{%
g\setminus }$ $\mathfrak{g}_{0}.$ Furthermore, since $[\mathfrak{g}_{r-1},%
\mathfrak{g]\subset g}_{r-2},$ we have $[x_{1},a]\in \mathfrak{g}_{r-2}.$
Now we argue in the same way replacing $a$ with $[x_{1},a]:$ If $[[x_{1},a],%
\mathfrak{g}]\subset \mathfrak{g}_{r-2},$ then $[x_{1},a]\in \mathfrak{g}%
_{r-1}$ by (8), which is not the case. Therefore, there exists some $%
x_{2}\in \mathfrak{g}$ such that $[x_{2},[x_{1},a]]\notin \mathfrak{g}%
_{r-2}. $ Clearly $x_{2}\in \mathfrak{g\setminus g}_{0}$ since $[\mathfrak{g}%
_{0},\mathfrak{g}_{r-2}]\subset \mathfrak{g}_{r-2}.$ Furthermore, $%
[x_{2},[x_{1},a]]\in \mathfrak{g}_{r-3}$ because $[\mathfrak{g}_{r-2},%
\mathfrak{g}]\subset \mathfrak{g}_{r-3}$ by (10). Iterating this process, we
obtain a sequence

\begin{equation}
a,\text{ }ad(x_{1})(a),\text{ }ad(x_{2},x_{1})(a),...,\text{ }%
ad(x_{r},...,x_{1})(a),\text{ }ad(x_{r+1},...,x_{1})(a)
\end{equation}%
with the following properties.

$1)$ $x_{1},...,x_{k+1}\in \mathfrak{g\setminus g}_{0}$

$2)$ $ad(x_{i},...,x_{1})(a)\in \mathfrak{g}_{r-i-1}\setminus \mathfrak{g}%
_{r-i},$ $1\leq i\leq r,$ where we set $\mathfrak{g}_{-1}=\mathfrak{g.}$

In particular, all the terms in (11) except the last term belong to $%
\mathfrak{g}_{0}.$ Now we observe that the term before the last one in (11)
belongs to $\mathfrak{g}_{(r-1)}=$ the $(r-1)$'th term in the lower central
series of $\mathfrak{g.}$ Indeed, let us recall that we set $\mathfrak{g}%
_{(0)}=\mathfrak{g}$ and define $\mathfrak{g}_{(k+1)}\overset{def}{=}[%
\mathfrak{g}_{(k)},\mathfrak{g].}$ Note the parantheses in $\mathfrak{g}%
_{(k)}$ to distinguish it from $\mathfrak{g}_{k}$ defined by (9). Also note
that $\mathfrak{g}_{(k)}$ depends only on $\mathfrak{g}_{(0)}=\mathfrak{g}$
whereas $\mathfrak{g}_{k}$ depends on the pair $(\mathfrak{g},\mathfrak{g}%
_{0})!$ Therefore we conclude

\begin{equation}
0\neq ad(x_{r},...,x_{1})(a)\in \mathfrak{g}_{0}\cap \mathfrak{g}_{(r-1)}
\end{equation}

Now (12) implies the following

\begin{proposition}
Let $\mathfrak{g}$ be a Lie algebra (always finite dimensional and over $%
\mathbb{R}$ or $\mathbb{C})$ and $0\neq \mathfrak{g}_{0}\subset \mathfrak{g}$
a subalgebra such that $(\mathfrak{g,g}_{0})$ is effective with order $r=o(%
\mathfrak{g,g}_{0}).$ Then $\mathfrak{g}_{0}\cap \mathfrak{g}_{(r-1)}\neq
\varnothing $ where $\mathfrak{g}_{(r-1)}=$ the $(r-1)$'th term in the lower
central series of $\mathfrak{g.}$
\end{proposition}

\begin{corollary}
Suppose $(\mathfrak{g,g}_{0})$ is effective with order $o(\mathfrak{g,g}%
_{0}) $ and $\mathfrak{g}$ is $k$-step nilpotent, i.e., $\mathfrak{g}%
_{(k-1)}\neq 0 $ but $\mathfrak{g}_{(k)}=0.$ Then $o(\mathfrak{g,g}_{0})\leq
k+1.$
\end{corollary}

For nilpotent $\mathfrak{g},$ we will denote the integer $k$ in Corollary 7
by $n(\mathfrak{g})$ and call it the nil-length of $\mathfrak{g.}$ According
to Corollary 7, a necessary condition to construct an effective $(\mathfrak{%
g,g}_{0})$ of high order with nilpotent $\mathfrak{g}$ is that the
nil-length of $\mathfrak{g}$ should be "large". Our purpose is now to give a
lower bound for order and search for a sufficient condition. Note that a
subalgebra $\mathfrak{h\subset }$ $\mathfrak{g}$ of a nilpotent $\mathfrak{g}
$ is also nilpotent and $n(\mathfrak{h})\leq n(\mathfrak{g}).\mathfrak{\ }$%
It is also well known (and easy to show) that if $\mathfrak{h\subsetneqq }$ $%
\mathfrak{g}$ with $\mathfrak{g}$ nilpotent, then $\mathfrak{h\subsetneqq }N(%
\mathfrak{h})=$ the normalizer of $\mathfrak{h.}$

\begin{lemma}
Let $\mathfrak{g}$ be nilpotent. If $\dim (\mathfrak{g})\geq 2$ and $%
\mathfrak{g}$ is nonabelian, then there exist proper subalgebras $\mathfrak{%
k\subsetneqq h\subsetneqq g}$ such that $\mathfrak{k}$ is normal in $%
\mathfrak{h,}$ $\mathfrak{h}$ is normal in $\mathfrak{g}$ but $\mathfrak{k}$
is not normal in $\mathfrak{g.}$
\end{lemma}

Let $(a)$ $=$ the linear span of $a$ and suppose that $(a)$ is a $1$%
-dimensional \textit{ideal} of $\mathfrak{g}$ for all $a\in \mathfrak{g.}$
Let $a,b\in \mathfrak{g.}$ If $a,b$ are linearly dependent, then clearly $%
[a,b]=0.$ Suppose they are linearly independent. Now $[a,b]=\lambda a=\mu b$
for some $\lambda ,\mu $ which implies $\lambda =\mu =0.$ Therefore $[a,b]=0$
for all $a,b\in \mathfrak{g,}$ contrary to our assumption. Thus there exists
some $a\in $ $\mathfrak{g}$ such that $(a)$ is not an ideal of $\mathfrak{g,}
$ i.e., $N(a)\subsetneqq \mathfrak{g.}$ Now letting $N_{0}(a)=(a),$ $%
N_{k+1}(a)=N(N_{k}(a)),$ we obtain the strictly increasing sequence $%
N_{1}(a)\subsetneqq N_{2}(a)\subsetneqq N_{3}(a)\subsetneqq ....$ If $k\geq
2 $ is the smallest integer with $N_{k}(a)=\mathfrak{g,}$ we set $\mathfrak{%
h=}N_{k-1}(a)$ and $\mathfrak{k=}N_{k-2}(a).$

In the next proposition we will need also the derived length of a solvable
Lie algebra $\mathfrak{g.}$ Recall that we set $\mathfrak{g}^{(0)}=$ $%
\mathfrak{g}$ and define inductively $\mathfrak{g}^{(k+1)}=[\mathfrak{g}%
^{(k)},\mathfrak{g}^{(k)}].$ Clearly $\mathfrak{g}^{(k)}\subset \mathfrak{g}%
_{(k)}.$ Therefore nilpotent implies solvable and the sol-length of $%
\mathfrak{g}=s(\mathfrak{g})\leq n(\mathfrak{g})=$ the nil-length of $%
\mathfrak{g}.$

\begin{proposition}
Let $\mathfrak{g}$ be nilpotent and suppose that the subalgebras $\mathfrak{%
k\subsetneqq h\subsetneqq g}$ satisfy the conclusion of Lemma 8. If $(%
\mathfrak{g,}\mathfrak{k)}$ is effective, then $o(\mathfrak{g,}\mathfrak{k}%
)\geq s(\mathfrak{k}).$
\end{proposition}

The first assumption of Proposition 9 is natural in view of Lemma 8.
However, there are two problems with Proposition 9. According to [S],
nilpotent Lie algebras with "sufficiently large" nil-length contain \textit{%
proper }subalgebras whose nil-length is also "large". Starting with such a
proper subalgebra, taking normalizers will not decrease the nil-lengths and
will eventually give two proper subalgebras satisfying the first assumption
of Proposition 9. The first problem is that these subalgebras may not have
also large sol-lengths (though this seems improbable). The second problem is
that we do know how to construct effective pairs using such subalgebras and
therefore simply assume their existence. Therefore, as it stands,
Proposition 9 is much too far from being a satisfactory answer to \textbf{Q'.%
}

The proof of Proposition 9 is easy. Clearly $[\mathfrak{k,g]\subset h}$ and $%
[\mathfrak{k,[k,g]]\subset k.}$ By the Jacobi identity, we conclude $[[%
\mathfrak{k,k],g]\subset k,}$ i.e., $\mathfrak{k}^{(1)}$ has the property $[%
\mathfrak{k}^{(1)}\mathfrak{,g]\subset k}.$ However, writing $\mathfrak{k}%
_{k}$ for $\mathfrak{g}_{k},$ $\mathfrak{k}_{1}$ is maximal with respect to
this property by (9). Hence we conclude $\mathfrak{k}^{(1)}\subset \mathfrak{%
k}_{1}.$ Now we continue inductively: Since $[\mathfrak{k}^{(1)}\mathfrak{%
,g]\subset k},$ we have $[\mathfrak{k}^{(1)}[\mathfrak{k}^{(1)}\mathfrak{%
,g]]\subset \lbrack k}^{(1)},\mathfrak{k]}$ $\subset \mathfrak{k}%
^{(1)}\subset \mathfrak{k}_{1}$ and therefore $[[\mathfrak{k}^{(1)},%
\mathfrak{k}^{(1)}]\mathfrak{,g]=[k}^{(2)},\mathfrak{g]\subset k}_{1}$ again
by Jacobi. Therefore $\mathfrak{k}^{(2)}\subset \mathfrak{k}_{2}$ by the
maximality of $\mathfrak{k}_{2}.$ Continuing in this way we get $\mathfrak{k}%
^{(k)}\subset \mathfrak{k}_{k}.$ Now $\mathfrak{k}^{(k)}$ stabilizes at zero
after $s(\mathfrak{k})$ steps and $\mathfrak{k}_{k}$ stabilizes at zero 
\textit{by assumption} which gives the required lower bound.

The above arguments give the impression that $o(\mathfrak{g,g}_{0})$ can be
large only if $\mathfrak{g}$ is nilpotent but this is by no means true.
Suppose the lower central series of a Lie algebra $\mathfrak{g}$ stabilizes
at the ideal $\mathfrak{h}$ with length $n(\mathfrak{g})\geq 0.$ For
instance $n(\mathfrak{g})=0$ if and only if $\mathfrak{g}^{(1)}=\mathfrak{g}%
_{(1)}=\mathfrak{g}.$ Such Lie algebras include semisimple ones and are
sometimes called perfect. Now $\mathfrak{g/h}$ is nilpotent and we can argue
with $\mathfrak{g/h}$ as above. At any rate, it is clear by now that for $(%
\mathfrak{g,g}_{0})$ to have a large order, $\mathfrak{g}$ and $\mathfrak{g}%
_{0}$ must have large nil-lengths and sol-lengths.

\bigskip

\bigskip

\bigskip

\textbf{References}

\bigskip

[GS] Guillemin, V., Sternberg, S., An algebraic model of transitive
differential geometry, Bull. Amer. Math. Soc. 70, 1964, 16-47

[O1] Ortacgil, E., An Alternative Approach to Lie Groups and Geometric
Structures, Oxford University Press, 2018

[O2] Orta\c{c}gil, E., On a simple proof of the Poincare Conjecture

[S] Safa, H., A bound for the nilpotency class of a Lie algebra, Iranian
Journal of Mathematical Sciences and Informatics Vol. 14, No. 2 (2019), pp
153-156

\bigskip

Erc\"{u}ment H. Orta\c{c}gil

ortacgile@gmail.com

\end{document}